\documentclass[11pt,letterpaper]{amsart}

\usepackage{latexsym}

\usepackage{amssymb}
\usepackage{amsthm}
\usepackage{amsmath}
\usepackage{amstext}
\usepackage{amscd}

\numberwithin{equation}{section}

\theoremstyle{plain}
\newtheorem{thm}{Theorem}[section]
\newtheorem{prop}[thm]{Proposition}
\newtheorem{cor}[thm]{Corollary}

\newtheorem{theorem*}{Theorem}[]

\theoremstyle{definition}
\newtheorem{defn}[thm]{Definition}

\newtheorem{example}[thm]{Example}

\theoremstyle{remark}
\newtheorem{rem}[thm]{Remark}

\newcommand{\C}{\mathbb{C}}

\newcommand{\R}{\mathbb{R}}

\newcommand{\Z}{\mathbb{Z}}

\newcommand{\Ko}{K_0(\mathcal V_{\R})}

\newcommand{\Ks}{K_0(\mathcal S)}
\newcommand{\chic}{\chi_c}
\newcommand{\en} {\varphi_n}
\newcommand{\bl}{\operatorname{Bl}}
\newcommand{\fund}[1]{\mu_{#1}}

\newcommand{\zar}[1]{\overline{#1}^Z}

\date{October 22, 2002}

\title{Virtual Betti numbers of real algebraic varieties}

\author{ Clint McCrory and Adam Parusi\'nski }

\address {D\'epartement de Math\'ematiques, Universit\'e d'Angers,
2 bd Lavoisier, 49045 Angers Cedex, France}

\email{parus@tonton.univ-angers.fr}

\address {Mathematics Department, University of Georgia, Athens GA
30602, USA}

\email{clint@math.uga.edu}

\begin{document}
 
\begin{abstract} The weak factorization theorem for birational maps
is used to prove that for all $i\geq 0$ the $i$th
mod 2 Betti number  of compact nonsingular real algebraic varieties
has a unique extension to a
\emph{virtual Betti number} $\beta_i$ defined for all real
algebraic varieties, such that if  
$Y$ is a
closed subvariety of
$X$, then 
$\beta_i(X) =
\beta_i(X\setminus Y) +  \beta_i(Y)$.
\end{abstract}

\maketitle


We define an invariant of the Grothendieck ring $\Ko$ of real
algebraic varieties, the \emph{virtual Poincar\'e polynomial}
$\beta(X,t)$, which is a ring homomorphism $\Ko\to\Z[t]$. For
$X$ nonsingular and compact,
$\beta(X,t)$ is the classical Poincar\'e polynomial for 
cohomology with $\Z_2$ coefficients. The
coefficients of the virtual Poincar\'e polynomial are the virtual
Betti numbers. By the weak factorization
theorem for birational morphisms (Abramovich \emph{et al.}
\cite{AKMW}), the existence of $\beta(X,t)$ follows from a simple
formula for the Betti numbers of the blowup of a compact
nonsingular variety along a closed nonsingular center.  The existence of the virtual Betti numbers for certain
real analytic spaces, including real algebraic varieties, has also been
announced by Totaro \cite{totaro}.

Kontsevich's motivic measure on the arc space of a complex algebraic
variety takes values in the completion of the localized Grothendieck
ring (\emph{cf.} \cite{looijenga}, \cite{denefloeser1},
\cite{craw}).
 Completion with respect to virtual dimension is possible because
complex varieties which are equivalent in the Grothendieck ring
have the same dimension. Using the virtual Poincar\' e polynomial we
prove that dimension is an invariant of the Grothendieck ring of
real algebraic varieties.
 It follows that motivic measures can be defined on arc spaces of
real varieties.

By a real algebraic variety we mean an algebraic variety in the
sense of Serre
\cite{serre} over the field $\R$ of real numbers. For our purposes
there is no loss of generality in restricting our attention to
affine real algebraic varieties; over $\R$ every quasiprojective
variety is affine. For background on real algebraic varieties we
refer the reader to \cite{BCR}.

\section{Generalized Euler characteristics}
\label{Euler}

\begin{defn}\label{groth} The \emph{Grothendieck group} of real
algebraic varieties is the
 abelian group generated by symbols $[X]$, where $X$ is a real
algebraic  variety, with the two relations
\begin{enumerate}
\item [(1)]
$[X] = [Y]$ if $X$ and $Y$ are isomorphic,
\item [(2)]
$[X] = [X \setminus Y] + [Y]$ if $Y$ is a closed subvariety of $X$.
\end{enumerate} The product of varieties induces a ring structure
\begin{enumerate}
\item [(3)]
$[X]\cdot [Y] = [X\times Y]$,
\end{enumerate} and the resulting ring, denoted by $\Ko$, is
called  the \emph{Grothendieck ring} of real algebraic varieties.
The class of a point is the unit of $\Ko$. The zero element of
$\Ko$ is the class of the empty set.
\end{defn} 

\begin{defn}\label{euler} A \emph{generalized (proper) Euler
characteristic} of real algebraic varieties, with values in a ring
$R$, is a ring homomorphism
$e:\Ko \to R$. In other words, for all varieties $X$ and $Y$,
\begin{enumerate}
\item [(1)]
$e(X) = e(Y)$ if $X$ and $Y$ are isomorphic,
\item [(2)]
$e(X) = e(X \setminus Y) + e(Y)$ if $Y$ is a closed subvariety of
$X$,
\item [(3)]
$e(X\times Y) = e(X)\cdot e(Y).$
\end{enumerate}
\end{defn}

The Euler characteristic with compact supports $\chic$ is an
example of such a homomorphism.  In fact the only generalized Euler
characteristics which are homeomorphism invariants are constant
multiples of
$\chic$ (\emph{cf.} \cite{quarez}).  The standard Euler
characteristic $\chi$ does not satisfy the additivity property
\ref{euler}(2) for  real algebraic varieties, and hence it is
\emph{not} a generalized Euler characteristic.

For example, let $X$ be the unit circle in the plane, and let $Y$ be
a  point of $X$. Then $\chi(X) = 0$, but
$\chi(X\setminus Y) = 1$ and $\chi(Y) = 1$. On the other hand,
$\chic(X) = 0$, 
$\chic(X\setminus Y) = -1$ and $\chic(Y) = 1$.

There are several ways to view the Euler characteristic with
compact supports. By definition $\chic(X)$ is the Euler
characteristic of the sheaf cohomology of
$X$ with compact supports, with coefficients in the constant sheaf
$\Z$. (Below we use coefficients $\Z_2$, the integers mod 2.) We
denote by $H^i_c(X)$ the
$i$th cohomology of $X$ with compact supports. Classical cohomology
has arbitrary closed supports. Property \ref{euler}(2) for $\chic$
follows from the long exact cohomology sequence of the pair
$(X,Y)$ (\cite{godement}, 4.10). If $X$ is triangulated then the
sheaf cohomology of $X$ with compact supports is isomorphic to the
simplicial cohomology of $X$ with compact supports.

For coefficients in a field the $i$th cohomology of $X$ with
compact supports is isomorphic to the $i$th Borel-Moore homology of
$X$ (\cite{borelmoore}, 3.3). Thus $\chic(X)$ is also equal to the
Euler characteristic of the Borel-Moore homology of $X$.
Borel-Moore homology has arbitrary closed supports. Classical
homology has compact supports. If $X$ is triangulated the
Borel-Moore homology of $X$ is isomorphic to the simplicial
homology of
$X$ with closed supports (\emph {i.e.} using possibly infinite
simplicial chains). If
$\bar X$ is a compact space containing $X$ such that $\bar X$ is
triangulated with
$\bar X\setminus X$ as a subcomplex, then the Borel-Moore homology
of $X$ is isomorphic to the relative simplicial homology of the pair
$(\bar X,\bar X \setminus X)$. Every real algebraic variety $X$
admits a compactification with such a triangulation (\emph{cf.}
\cite{lojasiewicz}, \cite{hironaka}, \cite{BCR}).

Let $X$ be a compact nonsingular real algebraic variety and let
$C$ be a closed nonsingular subvariety of $X$. Denote by $\bl _C X
\to X$ the  blow-up of $X$ with center $C$ and by $E$ the
exceptional divisor.  Then since $\bl_C X \setminus E$ is
isomorphic to $X\setminus C$, we have
\begin{equation}\label{blowup} [\bl_C X]- [E] = [X] - [C].
\end{equation} Bittner \cite{bittner} has proved using the weak
factorization theorem of Abramovich \emph{et al.}
\cite{AKMW} that one may replace
\ref{groth}(2) in the definition of the Grothendieck ring by
\eqref{blowup}.  The following result is equivalent to Bittner's,
and our proof is a reorganization of hers.

\begin{thm} \label{bittner} Let $X\mapsto e(X)$ be a function
defined on  compact nonsingular real algebraic varieties  with
values in a abelian group $G$, such that $e(X)=e(Y)$ if $X$ and
$Y$ are isomorphic.  Suppose that
$e(\emptyset) = 0$ and  for any blow-up $\bl _C X \to X$, where $X$
is a compact nonsingular real algebraic variety and $C\subset X$ a
closed nonsingular  subvariety, with exceptional divisor $E$,
\begin{equation*} e(\bl_C X) - e(E) = e(X) - e(C).
\end{equation*} Then $e$ extends uniquely to a group homomorphism
$\Ko \to G$.

If $e$ takes values in a ring $R$ and
$e(X\times Y) = e(X)\cdot e(Y)$ for all $X$, $Y$ compact nonsingular
 then this extension of $e$ is  a ring homomorphism $\Ko \to  R$.
\end{thm}

\begin{proof} First we note that if $X$ is the disjoint union of
the compact nonsingular varieties $Y$ and $Z$,
\begin{equation*} X = Y \sqcup Z,\end{equation*} then
 \begin{equation*} e(X) = e(Y) + e(Z).\end{equation*} For the
blow-up of $X$ with center $C=Y$ is $Z$, and the exceptional
divisor $E$ is empty (\emph {cf.} \cite{BCR}, p. 78 \emph {et
seq.}).

We prove the following statement by induction on $n$:

\noindent ($F_n$)
\emph{There is a unique map $X\mapsto \en(X)\in G$ defined on
varieties of dimension $\le n$  such that
\begin{enumerate}
\item[(a)] $\en(X)=\en(Y)$ if $X$ and
$Y$ are isomorphic,
\item [(b)]
$\en (X) = e(X)$ for $X$ compact and nonsingular,
\item [(c)]
$\en (X\setminus Y) = \en (X) - \en (Y)$ if $Y$ is a closed
subvariety of $X$. 
\end{enumerate}}

\smallskip The case $n=0$ is trivial.  Suppose ($F_{n-1}$) holds.
If $\dim X\le n-1$ then we put $\en (X) := \varphi_{n-1} (X)$.  If
$\dim X= n$ then $\en (X)$ is defined follows:
\begin{enumerate}
\item [(i)] Let $X$ be nonsingular and let $\bar X$ be a
nonsingular compactification of $X$.  Then
$$
\en (X) := e(\bar X) - \en (\bar X\setminus X),
$$ 
\item [(ii)] Let $X$ be singular and let $X= \bigsqcup S_i$ be a
stratification of
$X$. Then 
$$
\en (X) := \sum \en (S_i).
$$
\end{enumerate}

Here by a nonsingular compactification of the variety $X$ we mean a
compact nonsingular variety $\bar X$ which contains a subvariety
$U$ isomorphic to $X$ such that the  Zariski closure of $U$ in
$\bar X$ equals $\bar X$.

By a stratification $\mathcal S=\{S_i\}$ of
$X$ we mean a finite disjoint  union
$X= \bigsqcup S_i$ where each stratum  $S_i$ is a locally closed
nonsingular subvariety of $X$.  Denote by $\zar {S_i}$ the Zariski
closure of $S_i$ in $X$.  We also require
 that, for every stratum $S_i$, the set $\zar {S_i} \setminus S_i$
is a closed algebraic subvariety of $X$ which is a  union of strata
of $\mathcal S$.

First we show that the definition (i) of $\en(X)$ for $X$
nonsingular  is  independent of the choice of compactification. Let
$\bar X$ and $\bar X'$ be two nonsingular compactifications of
$X$. Then $\bar X$ and $\bar X'$ are birationally isomorphic.  By
the weak factorization  theorem there is a finite sequence of
blow-ups and blow-downs with  nonsingular centers joining $\bar X$
and $\bar X'$:
\begin{equation*}
\minCDarrowwidth 25pt
\begin{CD}
\bar X @<\varphi_0<< X_1 @>\varphi_1 >> X_2  @<<<\cdots @>>>
X_{k-1} @<\varphi_{k-1}<<    X_k @>\varphi_k >> \bar X'
\end{CD}
\end{equation*} Thus to prove  independence we may suppose that
$\bar X' = \bl _C \bar X$, where
$C$ is nonsingular and $C\cap X=\emptyset$. Let $E$ denote the
exceptional divisor.  Then
\begin{eqnarray*} & & e(\bar X')  - \en (\bar X' \setminus X) =
e(\bar X') - \en (E) -
\en (\bar X' \setminus X \setminus E)  \\ & & \quad = e(\bar X') -
\en (E) - 
\en (\bar X \setminus X \setminus C) = e(\bar X') - \en (E) +\en
(C) - 
\en (\bar X \setminus X) \\ & & \quad = (e(\bar X') - e (E) + e
(C))  - \en (\bar X\setminus X) =  e(\bar X)  - \en (\bar
X\setminus X) .  
\end{eqnarray*} 

Next we show that $\en (X\setminus Y) = \en (X) -\en (Y)$ for $X$
and $Y$ nonsingular. Let $\bar X$ be a nonsingular compactification
of $X$.   Let $D = \bar X \setminus X$. (It is not necessary to
assume that $D$ is a divisor in $\bar X$.)   First consider the
case 
$\zar{X\setminus Y} = X$, \emph{i.e.} $Y$ contains no irreducible
component of $X$.   Since $\dim Y < n$, we have
\begin{eqnarray*}
 \en (X\setminus Y) = e(\bar X) - \en (D\cup Y) = e(\bar X) - \en
(D) - 
\en (Y) = \en (X) -\en (Y) .
\end{eqnarray*} The next case is when $Y$ is the union of some
irreducible components of $X$.  Let $\bar Y$ denote the Zariski
closure of $Y$ in $\bar X$. Clearly $\bar Y$ is a union of
irreducible components of $\bar X$.  Let $W= \bar X \setminus \bar
Y$.
 Then
\begin{eqnarray*} & \en (X\setminus Y) & = e(W) - \en (D\cap W) \\
& & = e(\bar X) - e(\bar Y) - \en (D)  + \en (D\cap \bar Y) =
\en (X) - \en (Y) .
\end{eqnarray*} In the general case we denote by $X'$ the union of
the irreducible components of $X$ contained in $Y$ and by $X''$ the
union of the remaining irreducible components.  Let $Y'' = Y\cap
X''$.  Then by the above cases
\begin{eqnarray*} & \en (X\setminus Y) &= \en (X'' \setminus Y'') =
\en (X'') - \en (Y'') \\  & & = \en (X) - \en (X') - (\en (Y) - \en
(X'))  = 
\en (X) - \en (Y) .
\end{eqnarray*}

In particular, suppose $X$ is nonsingular, $\dim X\le n$, and let
$X= \bigsqcup S_i$ be a stratification of $X$.  Let $S_0$ be a
stratum of  smallest dimension.  Then $S_0$ is a closed nonsingular
variety of $X$ and hence
$$
\en (X\setminus S_0) = \en (X) - \en (S_0).
$$ Since $X\setminus S_0 = \bigsqcup_{S_i \ne S_0} S_i$ is again a
stratification, by induction on the number of strata we have
$$
\en (X) = \sum \en (S_i).
$$

Let $X$ be an arbitrary singular variety with $\dim X \le n$.
Consider two stratifications $X=\bigsqcup S_i$, $X=\bigsqcup S'_j$
of $X$. Suppose the second  stratification refines the first; that
is, each stratum $S_i$  is a union of strata $S'_j$. Then by the
additivity of
$\en$ for nonsingular varieties, we have
$\sum \en (S_i)= \sum \en (S'_j)$.   In general, given two
stratification  of $X$ there exists a third stratification refining
both of them.  This shows that the definition (ii) of
$\en (X)$  is independent of the stratification.

We now show the additivity (c) of $\en$. Suppose $\dim X \leq n$,
and let 
$Y$ be a subvariety of $X$.  There exists a  stratification of $X$
such that $Y$ is a union of strata.   So the additivity of $\en$
follows from the additivity of $\en$ for nonsingular  varieties.
This completes the inductive proof of ($F_n$).

Clearly the $\en$ defined above, $n\geq 0$, give a unique extension
of $e$ to a group homomorphism
$\varphi:\Ko \to G$. If $e$ is a ring homomorphism then it is easy
to see from (i) and (ii) above, by induction on dimension, that
$\varphi$ is also a ring homomorphism.
\end{proof} 


\section{The virtual Poincar\'e polynomial}
\label{poincare}

Let $X$ be a nonsingular compact real algebraic variety, and for
$i\geq 0$ let
$b_i(X)$ be the $i$th Betti number with $\Z_2$ coefficients,
\begin{eqnarray*} b_i(X) = \dim_{\Z_2} H^i(X;\Z_2) = \dim_{\Z_2}
H_i(X;\Z_2).
\end{eqnarray*}

\begin{prop} \label{bi} Let $\pi: \bl _C X \to X$ be the blow-up of
the nonsingular compact variety $X$ with nonsingular center
$C\subset X$.  Let $E$ denote the exceptional divisor. Then, for all
$i$,  
\begin{equation*} b_i (\bl_C X) - b_i (E) = b_i (X) - b_i (C).
\end{equation*}
\end{prop}

\begin{proof} Without loss of generality we may assume that $X$ has
pure dimension $n$ and
$\dim C <n$.  

Let $\tilde X= \bl _C X$.  Recall that $H_c^i(X;\Z_2)$ denotes the
$i$th mod 2 cohomology of $X$ with  compact supports.
 Since $\tilde X \setminus E \simeq X\setminus C$,  the diagram of
long exact sequences of closed embeddings
\begin{equation*}\minCDarrowwidth 1pt\begin{CD}
\cdots @>>> H^{i-1}(E;\Z_2) @>>> H^i_c (\tilde X\setminus E;\Z_2)
@>>>  H^i (\tilde X;\Z_2)
 @>>> H^i(E;\Z_2) @>>> \cdots\\ @. @AA\pi^*A @A\simeq A\pi^*A
@AA\pi^*A @AA\pi_*A @.\\
 \cdots @>>> H^{i-1}(C;\Z_2) @>>> H^i_c (X\setminus C;\Z_2) @>>>
H^i ( X;\Z_2)
 @>>> H^i(C;\Z_2) @>>> \cdots
\end{CD}\end{equation*}  induces (by a diagram chase) a long exact
sequence 
\begin{equation}\label{exact}
\cdots\to H^{i-1} (E;\Z_2) \to H^i (X;\Z_2)
\to H^i (C;\Z_2) \oplus  H^i (\tilde X;
\Z_2) \to H^i (E;\Z_2) \to\cdots  .
\end {equation}

Let $\fund {X}\in H_n( X; \Z_2)$ denote the fundamental class of
$X$.  We have $\pi_{*} (\fund {\tilde X})=
\fund {X}$. 
  Hence by  Poincar\'e duality 
$\pi^* : H^i(X; \Z_2) \to H^i(\tilde X; \Z_2)$
 is injective for all $i\geq 0$. (Let $a\in H^i(X;\Z_2)$. Now
$\pi_*(\pi^*a\frown\fund{\tilde X}) = a\frown\pi_*\fund{\tilde X}=
a\frown\fund X$. Thus $\pi^*a = 0$ implies
$a\frown\fund X = 0$, so $a = 0$ by Poincar\'e duality.) This shows
that \eqref{exact} splits  into short exact sequences
\begin{equation*}\label{shortexact} 0 \to H^{i} (X;\Z_2) \to H^i
(C;\Z_2) \oplus  H^i (\tilde X; \Z_2)
\to H^i (E;\Z_2) \to 0 ,
\end {equation*} so
\begin{equation*} b_i(C) + b_i(\tilde X) = b_i(X) + b_i(E),
\end{equation*} as desired.
\end{proof}

\begin{cor} For each nonnegative integer $i$ there exists a unique
group homomorphism
$\beta_i:  \Ko \to \Z$ such  that $\beta_i (X) = b_i(X)$ for $X$
compact nonsingular.

There exists a unique ring  homomorphism $\beta (\cdot,t): \Ko \to
\Z[t]$ such  that $\beta (X,t) = \sum_i b_i(X)t^i$ for $X$ compact
nonsingular.
\end{cor}

\begin{proof} By the K\"unneth formula, the classical Poincar\'e
polynomial $P(X,t) = \sum_i b_i(X)t^i$ is multiplicative for compact
nonsingular varieties.  Thus the corollary follows from Proposition
\ref{bi} and Theorem \ref{bittner}.
\end{proof}

\begin{defn} The integer $\beta_i(X)$ is the $i$th \emph{virtual
Betti number} of the real algebraic variety $X$, and the polynomial
$\beta(X,t)= \sum_i \beta_i(X)t^i$ is the \emph{virtual Poincar\'e
polynomial} of $X$.
\end{defn}

\begin{thm}\label{dimension} The virtual Poincar\'e polynomial
$\beta(X,t)$ is of degree $n=\dim X$, and $\beta_n(X) >0$.   In
particular, $[X]=[Y]$ implies  $\dim X = \dim Y$, and
$[X]\ne 0$ if $X\ne \emptyset$.  \end{thm}

\begin{proof} The proof of Theorem \ref{bittner}  gives us formulas
for $\beta_i(X)$ in terms of  a nonsingular compactification of
$X$ (\ref{bittner})(i), or in terms of a stratification of $X$
(\ref{bittner})(ii).  Using these formulas we establish the theorem
by induction on the dimension of $X$. 

If $\dim X = 0$ then $X$ is a finite set of points, so $\beta(X,t)
= a$, where $a$ is the number of points of $X$.

Suppose the theorem is true for varieties of dimension less than
$n$, and $X$ has dimension $n$. First suppose $X$ is nonsingular.
Let
$\bar X$ be a nonsingular compactification of $X$, and let $D=\bar X
\setminus X$. Then
\begin{equation*}
\beta(X,t) = \beta(\bar X, t) - \beta(D,t).
\end{equation*} Since $\dim D< n$, the polynomial $\beta(D,t)$ has
degree less than $n$ by inductive hypothesis. Since
$\bar X$ is compact and nonsingular, $\beta(\bar X, t)$ is the
Poincar\'e polynomial of $\bar X$, which is of degree $n$, with
$\beta_n(\bar X) =b_n(\bar X) > 0$. Thus $\beta(X,t)$ has degree
$n$ and
$\beta_n(X)> 0$.

Now if $X$ is an arbitrary variety of dimension $n$, let 
$X=\bigsqcup S_i$ be a stratification of $X$ by nonsingular
varieties
$S_i$. Then
\begin{equation*}
\beta(X,t) = \sum_i\beta(S_i,t).
\end{equation*} Since each stratum $S_i$ is nonsingular of
dimension at most $n$, we have that $\beta(S_i,t)$ has degree at
most $n$, and if  
$\beta(S_i,t)$ has degree $n$ then $\beta_n(S_i)>0$. So
$\beta(X,t)$ has degree $n$ and $\beta_n(X)>0$.
\end{proof}

\begin{rem} In \cite{quarez} Quarez studies the
Grothendieck ring $\Ks$  of  semialgebraic sets, which is generated
by homeomorphism classes of semialgebraic sets with sum relation
\ref{groth}(2) and  product given by
\ref{groth}(3).  He observes that  $[X]=[Y]$ in $\Ks$ if and only
if $\chic(X) =\chic(Y)$.  In particular, the class  of a non-empty
semialgebraic set can be zero in $\Ks$, and two  semialgebraic sets
of different dimensions may represent the same  class in $\Ks$. 
This makes the construction of motivic measures---more precisely, 
completion with respect to virtual dimension---impossible for
$\Ks$.  
\end{rem}

The generalized Euler characteristics
\begin{eqnarray*} &\beta(X, -1)  = \sum_i (-1)^i \beta_i(X) \\
&\chic (X)  = \sum_i (-1)^i \dim H^i_c(X;\Z_2)
\end{eqnarray*} are equal for $X$ compact and nonsingular---they
both equal the Euler characteristic $\chi (X)$. It follows from
Theorem \ref{bittner} that
$\beta(X,-1)=\chic(X)$ for all real algebraic varieties $X$. But in
general $\beta_i(X) \ne \dim H^i_c(X;\Z_2)$. In fact $\beta_i(X)$
can be negative for $i < \dim X$, and the virtual Betti numbers
$\beta_i(X)$ are not topological invariants.

\begin{example} Consider the union of two intersecting ellipses,
\begin{equation*} X = \{(x,y)\ |\ (2x^2 + y^2 - 1)(x^2 + 2y^2-1) =
0\}.
\end{equation*} Let $X_1$ and $X_2$ be the two irreducible
components of $X$.  Then 
\begin{equation*}
\beta_0(X) = \beta_0(X_1) + \beta_0(X_2) - \beta_0(X_1\cap X_2) =
-2.
\end{equation*}
\end{example}

\begin{example} For the ``figure eight'' curve
\begin{equation*} X = \{(x,y)\ |\ y^2 = x^2 - x^4\},
\end{equation*} the proper transform of $X$ under the blowup of the
plane at the origin is homeomorphic to a circle, and the preimage
of the singular point of $X$ is two points. It follows that
$\beta_1(X) = 1$.  Now consider a second ``figure eight'' curve
\begin{equation*} Y = \{(x,y)\ |\ ((x+1)^2+y^2-1)((x-1)^2+y^2-1) =
0\}.
\end{equation*} Since $Y$ is the union of two circles tangent at
the origin, it follows that $\beta_1(Y) = 2$. But $Y$ is
homeomorphic to $X$.
\end{example}

\begin{example} Let $S^n$ be the unit $n$-sphere in $\R^{n+1}$. 
From the inclusion $\R^{n+1}\subset\R^{n+2}$ we have $S^n\subset
S^{n+1}$. From the decomposition $S^{n+1} = (S^{n+1}\setminus
S^n)\sqcup S^n$ we see that $\beta_n(S^{n+1}\setminus S^n) = -1$.
Now $S^{n+1}\setminus S^n$ is homeomorphic to
$\R^{n+1}\sqcup\R^{n+1}$, but $\beta_n(\R^{n+1}\sqcup\R^{n+1}) =
0$. For $S^{n+1}$ is the Alexandroff compactification of $\R^{n+1}$
(\emph {cf.} \cite{BCR}, p. 76), so $\beta_n(\R^{n+1})= 0$.
\end{example}


\section{Complex varieties}
\label{complex}

If we apply the complex versions of Theorem \ref{bittner} and
Proposition \ref{bi} to the Poincar\'e polynomial
\begin{equation*} P_{\C}(X,t)= \sum_i \dim_{\C} H^i(X;\C)t^i,
\end{equation*} we obtain  a generalized Euler characteristic of
complex algebraic varieties. This invariant is related to Deligne's
mixed Hodge theory in the following way.
   
For a complex algebraic variety $X$, the \emph{$E$-polynomial} (or
\emph{Hodge number characteristic}) is given by
\begin{equation*} E(X,u,v) = \sum_{i,p,q}(-1)^ih_{p,q}^iu^pv^q,
\end{equation*}
  where $h_{p,q}^i$ is  the dimension of the $(p,q)$-component of
the mixed Hodge structure on
$H^i_c(X;\C)$. The $E$-polynomial is a generalized Euler
characteristic (see \cite{danilovkhovanskii}, \cite{durfee},
\cite{looijenga}, \cite{craw}). Therefore so is the \emph{weight
characteristic} 
\begin{equation}\label{weight} E(X,t,t) = \sum_{i,j}(-1)^iw_j^i t^j.
\end{equation} The coefficients of the weight characteristic are
given by
\begin{equation*} w_j^i(X)=\sum_{p+q=j}h_{p,q}^i(X)=\dim_\C
W_j^i(X)/W_{j-1}^i(X),
\end{equation*} where 
\begin{equation}\label{filtration} 0\subset W^i_0(X) \subset
W^i_1(X)
\subset \cdots \subset W^i_i(X) = H^i_c(X;\C)
\end{equation} is the \emph{weight filtration} of cohomology with
compact supports
\cite{deligne}.

 If $X$ is compact and nonsingular then $W^i_{i-1} = 0$, so
\begin{equation*} w_j^i(X) =
\begin{cases} b_i^\C(X)\ \ i=j\\
\ \ 0\ \ \ \ \   i\neq j,
\end{cases}
\end{equation*} 
where $b_i^\C(X)= \dim_\C H^i(X;\C)$,
and thus the weight characteristic is the Poincar\'e
polynomial $P_\C(X,-t)$. In other words, for complex algebraic
varieties the virtual Poincar\'e polynomial is the weight
characteristic evaluated at
$-t$, and so the $j$th virtual Betti number is $(-1)^j$ times
the \emph{weight $j$ Euler characteristic},
\begin{equation*}
\beta_j^{\C}(X)=(-1)^j\sum_i(-1)^iw^i_j(X).
\end{equation*} The weight $j$ Euler characteristic has been
studied by Durfee \cite{durfee}. The virtual Poincar\'e polynomial
has been used by Fulton (\cite{fulton}, p. 92) to compute Betti
numbers of toric varieties.

\begin{rem} Applying the complex versions of Theorem \ref{bittner}
and Proposition \ref{bi} to the Poincar\'e polynomial with
coefficients in $\Z_p$, $p$ prime,  we get further examples of
generalized Euler characteristics of complex algebraic varieties.
The existence of these invariants also follows from the work of
Gillet and Soul\'e \cite{gilletsoule}.
\end{rem}

\begin{rem} One could study the Grothendieck ring of real algebraic
varieties
 by means of complexification, that is by considering the
isomorphism classes of pairs
$(X_{\C},\tau)$, where $X_{\C}$ is a complex algebraic variety
 with complex conjugation $\tau$.  The real algebraic variety
$X_\R$ associated to
$(X_{\C},\tau)$ is the fix point set of $\tau$. We do not know how
the generalized Euler characteristics of
$X_\C$ and of $X_\R$ are related, except for the observation that
\begin{equation*}
\chic (X_\C) \equiv \chic (X_\R) \mod 2.
\end{equation*}
 Note that a given real algebraic variety admits many different
realizations as the fixed point set of conjugation on a complex
variety.
\end{rem}


\medskip
\section{Weight filtration for real varieties}
\label{weight}

It is natural to ask whether the virtual Betti numbers of real
algebraic varieties are associated to a weight filtration on mod 2
cohomology with compact supports. Totaro \cite{totaro} has announced
that this is the case. The following examples show that a
weight filtration for real varieties cannot have properties as
strong as those of the weight filtration for complex varieties.

We analyze two examples. Example \ref{example1}, a curve in $\R^3$,
shows that a real weight filtration cannot have both the strict
naturality property and the resolution of singularities property
enjoyed by the complex weight filtration. Example \ref{example2}, a
surface in $\R^3$, shows
that there is no natural real weight filtration such that the
virtual Betti numbers are the weighted Euler characteristics.

\begin{defn} A \emph{real weight filtration} $W$ assigns to every
real algebraic variety $X$, and to every $i\geq 0$, a filtration of
the $i$th cohomology of $X$ with compact supports and $\Z_2$
coefficients, of the form
\begin{equation*}
0\subset W^i_0(X) \subset W^i_1(X)\subset\cdots\subset
W^i_i(X)=H^i_c(X;\Z_2).
\end{equation*}
If $W$ is a real weight filtration, for all $i$, $j$ we let
\begin{equation*}
w^i_j(X) = \dim_{\Z_2}W^i_j(X)/W^i_{j-1}(X).
\end{equation*}
\end{defn}
We will consider the following properties of a real weight
filtration. The classical complex weight filtration
considered in the previous section has all of these properties.
\begin{enumerate}
\item [(1)]
\emph{Naturality}. If $f:X\to Y$ is an algebraic morphism, then for
all
$i$ and $j$,
\begin{equation*}
f^*W^i_j(Y)\subset W^i_j(X).
\end{equation*}
\item[(2)]
\emph{Strict naturality}. If $f:X\to Y$ is an algebraic morphism,
then for all $i$ and $j$,
\begin{equation*}
f^*W^i_j(Y)= W^i_j(X)\cap \operatorname {Im} f^*.
\end{equation*}
\item[(3)]
\emph{Manifold}. If $X$ is compact and nonsingular, then for all
$i$
\begin{equation*}
W^i_{i-1}(X) = 0.
\end{equation*}
\item[(4)]
\emph{Resolution}. If $X$ is compact and $p:\tilde X\to X$ is a
resolution of singularities, then for all $i$
\begin{equation*}
W^i_{i-1}(X) = \operatorname{Ker} p^*.
\end{equation*} 
\item[(5)]
\emph{Virtual Betti}. For all $X$ and all $i$, $j\geq 0$, the
virtual Betti number $\beta_i$ is given by
\begin{equation*}
\beta_j(X)= (-1)^j\sum_i(-1)^iw^i_j(X).
\end{equation*}
\item[(6)]
\emph{Mayer-Vietoris}. If $A$ and $B$ are closed subvarieties of $X$
with
$A\cup B=X$ and $A\cap B = C$, then for all $j$ 
the Mayer-Vietoris cohomology sequence restricts to an exact
sequence
\begin{equation*}
\cdots\to W^i_j(X)\to W^i_j(A) \oplus W^i_j(B)\to W^i_j(C)\to
W^{i+1}_j(X)\to\cdots
\end{equation*}
\item[(7)]
\emph{Pair}. If $A$ is a closed subvariety of $X$, then 
for all $j$ 
the exact cohomology sequence of the pair restricts to an exact
sequence
\begin{equation*}
\cdots\to W^i_j(X\setminus A)\to W^i_j(X) \to
W^i_j(A)\to W^{i+1}_j(X\setminus A)\to\cdots
\end{equation*}
\end{enumerate}

Clearly condition (2) implies (1), (4) implies (3), (7)
and (3) imply (5). Since the virtual Betti numbers of a
compact nonsingular variety equal the classical Betti numbers, the
virtual Betti condition (5) implies the manifold condition (3).
\vskip.2in

In the rest of this section we will write $H^i(X) = H^i_c(X;\Z_2)$.
All the varieties we will consider are compact.

\begin{example}\label{example1}
In $\R^3$ with coordinates $(x,y,z)$, let $X$ be the intersection of
the circular cylinder
$x^2+y^2 = 1$ with the union of the two parabolic cylinders $x =
z^2$ and $x = -z^2$,
\begin{equation*}
X = \{(x,y,z)\ |\ x^2+y^2 = 1,\ (x-z^2)(x+z^2)=0 \}.
\end{equation*}
Thus $X$ is topologically the union of two circles which are tangent
at two points. Let $\tilde X$ be the disjoint union of the two
irreducible components of $X$, and let $p:\tilde X\to X$ be the
resolution of singularities given by the inclusions of the
components. Let
$C$ be the unit circle in the
$(x,y)$ plane, and let $q:X\to C$ be the projection. Now
$\dim H^1(X)=3$, and the sequence
\begin{equation*}
\minCDarrowwidth 1pt\begin{CD}
0 @>>> H^1(C)@>q^*>> H^1(X)@>p^*>> H^1(\tilde X)@>>>  0
\end{CD}
\end{equation*}
is exact. A real weight filtration on $H^1(X)$ has the form
\begin{equation*}
0\subset W^1_0(X)\subset W^1_1(X) = H^1(X).
\end{equation*}
If this filtration has the resolution property, then $W^1_0(X)=
\operatorname {Ker} p^* = \operatorname {Im} q^*$. Thus strict
naturality of the weight filtration implies that
$q^*W^1_0(C)=W^1_0(X)$. But the manifold property implies that
$W^1_0(C)= 0$, which is a contradiction.

Therefore a real weight filtration cannot satisfy both strict
naturality (2) and resolution (4). This example can also be used to
show that a real weight filtration cannot satisfy all three
of the following conditions: strict naturality, the manifold
condition, and the Mayer-Vietoris condition.
\end{example}

To analyze the next example we first restate the virtual Betti
condition (5) for surfaces.
Let $X$ be a 2-dimensional (compact) real algebraic variety. 
Consider a real weight filtration of $H^*(X)$,
\begin{eqnarray*}
& 0 \subset W^0_0(X) = H^0(X) \\
& 0 \subset W^1_0(X)\subset W^1_1(X) = H^1(X) \\
& 0 \subset W^2_0(X)\subset  W^2_1(X) \subset  W^2_2(X) = H^2(X).
\end{eqnarray*}
For $i\geq 0$ let $b_i$ be the $i$th mod 2 Betti
number of $X$, and let
$\beta_i$ be the
$i$th virtual Betti number of $X$.
The non-negative integers $w^i_j = \dim W^i_j(X)/W^i_{j-1}(X)$
satisfy the equations
\begin{eqnarray}
\nonumber &  w^0_0 = b_0 \\
\nonumber&   w^1_0  + w^1_1 = b_1\\
\label{wij} &    w^2_0  + w^2_1 + w^2_2 = b_2\\
\nonumber &   w^0_0  - w^1_0 + w^2_0 = \beta_0\\
\nonumber &   w^1_1  - w^2_1 =  \beta_1\\
\nonumber &    w^2_2 = \beta_2
\end{eqnarray}
Given $b_0$,  $b_1$, $b_2$, $\beta_0$,
$\beta_1$, $\beta_2$, we wish to determine the possible values of
$w^i_j$ which satisfy the 
equations (\ref{wij}). If we arrange the $w^i_j$ in an array
\vskip.2in
\begin{center}
\begin{tabular}{|c|c|c|}
	\hline
$w^2_2$ &  &  \\
	\hline 
 $w^1_1$ & $w^2_1$   & \\
	\hline
 $w^0_0$ & $w^1_0$ & $w^2_0$  \\
	\hline
\end{tabular}
\end{center}
\vskip.2in
\noindent
then
(\ref{wij}) says that the sum of the $i$th diagonal is $b_i$ and the
alternating sum of the $j$th row is
$\beta_j$ (\emph{cf.} the Mayer-Vietoris spectral sequence computation
below).

\begin{example}\label{example2}
Let $X=X_1\cup X_2\cup X_3$ be the divisor with normal
crossings in 3-space with smooth components $X_1$, $X_2$, $X_3$
defined as follows.

Let $X_1$ be the sphere of radius
$\sqrt 2$ with center $(1,0,0)$, 
\begin{equation*}
X_1=\{(x,y,z)\ |\ (x-1)^2 + y^2 + z^2 = 2\}.
\end{equation*}

 Let $X_2$ be the
sphere of radius $\sqrt 2$ with center $(-1,0,0)$,
\begin{equation*}
X_2=\{(x,y,z)\ |\ (x+1)^2 + y^2 + z^2 = 2\}.
\end{equation*}
Then $X_1\cap X_2$ is the circle of radius 1 in the $(y,z)$ plane
with center at the origin.

Let $X_3$ be the torus of revolution which is a tube of radius
$\sqrt 2/2$ with core the circle of radius 4 in the
$(y,z)$ plane with center $(0,4,0)$, 
\begin{equation*}
X_3=\{(x,y,z)\ |\ x^2 + (\sqrt{(y-4)^2+z^2} -4)^2 = 1/2\}.
\end{equation*}
($X_3$ is defined by a polynomial equation of degree 4.) 

The torus
$X_3$ intersects each of the spheres
$X_1$ and $X_2$ transversely along a simple closed curve. The
torus $X_3$ meets the circle $X_1\cap X_2$ transversely at four
points $p$, $q$, $r$, $s$. The locus $\Sigma$ of singularities of
$X$ is the union of the three loops $X_1\cap X_2$, $X_1\cap X_3$,
$X_2\cap X_3$, any two of which meet at the four points $p$, $q$,
$r$, $s$. The natural stratification of
$X$ consists of the nonsingular points $X\setminus \Sigma$, the
double points $\Sigma\setminus \{p,q,r,s\}$, and the triple points
$\{p,q,r,s\}$. The stratum $X\setminus \Sigma$ has 17 connected
components, 6 on each sphere and 5 on the torus. All but one of
these components is a topological 2-cell. The remaining component is
a topological annulus. The stratum $\Sigma \setminus \{p,q,r,s\}$
has 12 components, each of which is a 1-cell. So $X$ has a cell
structure with 4 0-cells, 13 1-cells, and 17 2-cells, with Euler
characteristic 8.

A computation with Mayer-Vietoris sequences (equivalent to the
Mayer-Vietoris spectral sequence computation below) gives
that $b_0 = 1$, $b_1 = 1$, and $b_2 = 8$.  The first homology
group of $X$ is generated by a circle on the torus $X_3$ which is
parallel to the core. The second homology group of
$X$ is generated by the boundaries of the 8 bounded connected
components of $\R^3\setminus X$. 

Now we compute the virtual Betti numbers of $X$. 
For $j=0, 1, 2$ we have
\begin{eqnarray*}
\lefteqn{\beta_j(X) =\beta_j(X_1)+\beta_j(X_2)+\beta_j(X_3)} \\
& &- \beta_j(X_1\cap X_2)-\beta_j(X_1\cap X_3) - \beta_j(X_2\cap
X_3)\\ & &+ \beta_j(X_1\cap X_2\cap X_3).
\end{eqnarray*}
Since the $X_i\cap X_j$ are
homeomorphic to circles and $X_1\cap X_2\cap X_3$ is four points,
we have
$\beta_0 = 4$, $\beta_1 = -1$, and $\beta_2 = 3$.

Let $W$ be a real weight filtration on the cohomology of $X$, and
suppose that $W$ is natural and satisfies the virtual Betti
condition. Since $b_0=1$, we have
$w^0_0=1$. Since $\beta_2 =3$ we obtain that $w^2_2 = 3$.

Since $b_1=1$, we have $w^1_0 + w^1_1 = 1$. Suppose that $w^1_0 =
1$. Then
$H^1(X) = W^1_0(X).$ But the inclusion $f:X_3\to X$ induces an
injection $f^*: H^1(X)\to H^1(X_3)$. Since $X_3$ is compact and
nonsingular, it follows from the virtual Betti condition that
$W^1_0(X_3) = 0$. This contradicts the naturality statement
$f^*W^1_0(X)\subset W^1_0(X_3)$.

The remaining option is that $w^1_1= 1$, which implies by
(\ref{wij}) that $w^2_1=2$ and $w^2_0=3$.

The subvariety $X_1\cup X_2$ of $X$ has $b_0=1$, $b_1=0$, $b_2=3$,
$\beta_0=1$, $\beta_1=-1$, $\beta_2=2$. By (\ref{wij}) we have
$w^2_0(X_1\cup X_2)=0$, $w^2_1(X_1\cup X_2)=1$, and $w^2_2(X_1\cup
X_2)=2$.

The subvariety $X_1\cup X_3$ has $b_0=1$,
$b_1=2$, $b_2=3$,
$\beta_0=1$, $\beta_1=1$, $\beta_2=2$. Now $w^1_0(X_3)=0$,
$w^1_1(X_3)=2$, and the restriction $H^1(X_1\cup X_3)\to H^1(X_3)$
is an isomorphism. Therefore, by naturality of $W$, $w^1_0(X_1\cup
X_3) = 0$. By (\ref{wij}) we conclude that $w^2_0(X_1\cup X_3)=0$,
$w^2_1(X_1\cup X_3)=1$, and $w^2_2(X_1\cup X_3)=2$. The same
computation applies to $X_2\cup X_3$.

There exists a basis $\{a_1,a_2,a_3,a_{12},a_{13},a_{23},b,c\}$ of
$H^2(X)$ such that each $a_i$ restricts to the generator of
$H^2(X_i)$, and each
$a_{ij}$ restricts to a generator of $W^2_1(X_i\cup X_j)$. (The class
dual to $a_i$ in homology is represented by the surface $X_i$, and
the class dual to $a_{ij}$ in homology is represented by the
boundary of the intersection of the region inside $X_i$ and
the region inside $X_j$.) Let $f_i:X_i\to X$ and $f_{ij}:X_i\cup
X_j\to X$, $1\leq i < j\leq 3$, be the inclusions. By naturality of
$W$, for
$k=0,1,2,$
\begin{equation*}
f_i^*W^2_k(X)\subset W^2_k(X_i),\ \ \ f_{ij}^*W^2_k(X)\subset
W^2_k(X_i\cup X_j).
\end{equation*}
It follows that the images of $a_1$, $a_2$, $a_3$ form a basis for
$W^2_2(X)/W^2_1(X)$, and the images of 
$a_{12}$, $a_{13}$, $a_{23}$ are independent in $W^2_1(X)/W^2_0(X)$.
Therefore $w^2_1\geq 3$, which is a contradiction.

Therefore a real weight filtration cannot satisfy both
naturality (1) and the virtual Betti condition (5). This example can
also be used to show that a real weight filtration cannot satisfy
both the Mayer-Vietoris condition and the virtual Betti condition.
\end{example}

It is instructive to compute for Example \ref{example2} the
Mayer-Vietoris spectral sequence which expresses the cohomology of
$X$ in terms of the cohomology of its irreducible components $X_1$,
$X_2$, $X_3$.  Let $X^0$ be the disjoint union of $X_1$,
$X_2$, and $X_3$. Let $X^1$ be the disjoint union of $X_1\cap X_2$,
$X_1\cap X_3$, and $X_2\cap X_3$. Let $X^2 = X_1\cap X_2\cap
X_3$. 

The Mayer-Vietoris spectral sequence is the spectral sequence of the
double complex $C^q(X^p)$, where $C^q$ is the group of
singular cochains (mod 2).  It converges to the
cohomology of $X$ (\emph{cf.} \cite{godement}, ch. 2, \S 5).
We have  $E_1^{pq} = H^q(X^p)$.
Since each
$X_i\cap X_j$ is homeomorphic to a circle and $X_1\cap X_2\cap
X_3$ is four points, we obtain the $E_1$ term:
\vskip.2in
\begin{center}
\begin{tabular}{|c|c|c|}
	\hline
$(\Z_2)^3$ &  &  \\
	\hline 
 $(\Z_2)^2$ & $(\Z_2)^3$   & \\
	\hline
 $(\Z_2)^3$ & $(\Z_2)^3$ & $(\Z_2)^4$  \\
	\hline
\end{tabular}
\end{center}
\vskip.2in
\noindent
Computing the differential
$d_1:E_1^{pq}\to E_1^{p+1,q}$ we find $E_2$:
\vskip.2in
\begin{center}
\begin{tabular}{|c|c|c|}
	\hline
$(\Z_2)^3$ &  &  \\
	\hline 
 $(\Z_2)^2$ & $(\Z_2)^3$   & \\
	\hline
 $\Z_2$ & $0$ & $(\Z_2)^3$  \\
	\hline
\end{tabular}
\end{center}
\vskip.2in
\noindent
The only possible non-zero differential $d_2^{0,1}:E^{0,1}_2\to
E^{2,0}_2$ has rank 1, and hence we obtain
$E_3$:
\vskip.2in
\begin{center}
\begin{tabular}{|c|c|c|}
	\hline
$(\Z_2)^3$ &  &  \\
	\hline 
 $\Z_2$ & $(\Z_2)^3$   & \\
	\hline
 $\Z_2$ & $0$ & $(\Z_2)^2$  \\
	\hline
\end{tabular}
\end{center}
\vskip.2in
\noindent
Finally, $d_r=0$ for $r\geq 3$, so $E_3=E_\infty$. This computation
gives that
$b_0= 1$,
$b_1=1$, and $b_2=8$, as claimed above.

The filtration $W$ of the cohomology of $X$
corresponding to the Mayer-Vietoris spectral sequence has
$w^i_j=\dim E^{i-j,j}_\infty$. The alternating sum
of the ranks of the entries of the $j$th row of $E_1$ is the
virtual Betti number $\beta_j$. Because the
$j$th row of $E_2$ is the homology of the $j$th row of $E_1$, the
alternating sum of the ranks of the entries of the $j$th row of
$E_2$ is also
$\beta_j$. (This holds for the Mayer-Vietoris
spectral sequence of any divisor with normal crossings.) Since the
differential $d_2$ is non-zero, the terms $E_r$ for
$r>2$ no longer have this property. Thus the filtration
associated to the Mayer-Vietoris spectral sequence does not
satisfy the virtual Betti condition.


\end{document}